\documentclass[10pt]{amsart}
\usepackage{latexsym}
\usepackage{amsmath}
\usepackage{amssymb}
\usepackage{mathrsfs}
\usepackage{graphicx}

\input epsf

\begin{document}

\title
{MATH-SELFIE}
\author{S.~S. KUTATELADZE}
\address{
Sobolev Institute of Mathematics\newline
\indent 4 Koptyug Avenue\newline
\indent Novosibirsk, 630090, Russia
}
\email{
sskut@math.nsc.ru
}
\begin{abstract}
This is a write up on some sections of
convex geometry, functional analysis,
optimization, and nonstandard models that
attract the author.
\end{abstract}

\maketitle

Mathematics is the logic of natural sciences, the unique
science of the provable forms of  reasoning
quantitatively and qualitatively.

Functional analysis had emerged at the junctions of
geometry, algebra and the classical calculus,
while turning rather rapidly into the natural language
of many traditional areas of continuous mathematics and
approximate methods of analysis. Also, it  has brought about
the principally new
technologies of theoretical physics and social sciences
(primarily, economics and control).

Of  most interest for the author are
some contiguous sections of the constituents of functional analysis and
mathematical logic that are promising in search for modernization of the theoretical
techniques of socializing the problems with many solutions
which involves the recent ideas of model theory.

The traditions of functional analysis
were implanted  in Siberia by S.~L. Sobolev and L.~V. Kantorovich.
Their thesis of the unity of functional analysis and applied mathematics
was, is, and should be the branding mark of the Russian mathematical school.
This is the author's deep belief.

The main areas dwelt upon below
are functional analysis, nonstandard methods of analysis, convex geometry, and optimization. These are listed according to their significance.
Each of them remains in the sphere of the author's interests from the moment
of the first appearance, but the time spent and the efforts allotted have been
changing every now and then.
Here these areas are addressed in chronological order.

\section{Optimal Location of Convex Bodies}

Using the ideas of linear programming invented by L.~V. Kantorovich,
it turned out possible to distinguish some classes of extremal problems
of optimal location of convex surface that  could not
be treated by the classical methods in principle.
The decisive step forward was to address such a problem by the
standard approach of programming which consists in transition  to
the dual problem. The latter turns out solvable by the technique of mixed
volumes,  abstraction of the duality ideas of H.~Minkowski, and
modification of one construction in measure theory that
belongs to Yu.~G. Reshetnyak.\footnote{Cp.~\cite{Reshetnyak}.}
The revealed descriptions of new classes of inequalities
over convex surfaces in combination with the technique of surface area measures
by A.~D. Alexandrov\footnote{Cp.~\cite{Alexandrov}.} had led to
 reducing to linear programs the isoperimetric-type
problems with however many constraints; i.e., the  problems
that fall beyond the possibility of symmetrizations.
In fact, the extensive class was discovered of the problems whose solutions can be written down
explicitly by translating the problems into convex programs in appropriate function
spaces.\footnote{Part of these results are presented in the survey
article~\cite{Duality}. The conception of $H$-convexity from this paper
is considered now as definitive in the numerous studies in generalized convexity
and search for schemes of global optimization
(in particular, cp.~\cite{Singer}).}

   The most visual and essential progress
is connected with studying some abstractions of the Urysohn problem
of maximizing the volume of a convex surface given the integral of
the breadth of the  surface. By the classical result of P.~S. Urysohn
which was published in the year of his death---1924,\footnote{Cp.~\cite{Urysohn}.}
this is a ball as follows from the suitable symmetry argument.
In the 1970s  the functional-analytical approach was illustrated
with the example of the {\it internal Urysohn problem}: Granted the integral breadth,
maximize the volume of a convex surface that lies within an a priori given
convex body, e.g., a simplex in~${\mathbb R}^N$.
The principal new obstacle in the problem is that
no symmetry argument is applicable in analogous internal or external problems.
It turned out that we may solve the problem in some generalized sense---``modulo''
the celebrated Alexandrov Theorem on reconstruction of a convex surface from its surface
area measure. For the Urysohn problem within a polyhedron
the solution will be given by the Lebesgue measure on the unit sphere
with extra point loads at the outer normals to the facets of the polyhedron.
The internal isoperimetric problem falls beyond the general scheme even
within a tetrahedron.

Considering the case of $N=3$ in 1995,
A.~V. Pogorelov found in one of his last
papers\footnote{Cp.~\cite{Pogorelov}.} the shape of the ``soup bubble''
within a tetrahedron in the same generalized sense---this happens to be
the vector sum of a ball and the solution of the internal Urysohn problem.
In the recent years quite a few papers has been written about
the double bubbles. These studies are also close to the
above ideas.

\section{Ordered Vector Spaces}

   Of most importance in this area of functional analysis are
the problems stemming from the Kantorovich heuristic principle.

In his first paper of 1935 on the brand-new topic L.~V. Kantorovich
wrote:\footnote{Cp.~\cite{LVK}.}

\smallskip
\hangindent 20pt{\small\quad\quad
In this note, I define a~new type of space that I call
a~semiordered linear space. The introduction of such a~space
allows us to
study linear operations of one abstract class
(those with values in these spaces) in the same way
as linear functionals.}
\medskip

It is worth noting that his definition of  semiordered
linear space contains the axiom of Dedekind completeness which
was denoted by $I_6$.  L.~V. Kantorovich demonstrated the role of
$K$-spaces by widening  the scope of the  Hahn--Banach Theorem.
The heuristic principle turned out applicable to this fundamental
Dominated Extension Theorem; i.e., we may abstract the Hahn--Banach Theorem on
substituting the elements of an arbitrary  $K$-space
for reals  and replacing linear functionals
with operators acting into the space.

The Kantorovich heuristic principle has found compelling justifications in his
own research as well as in the articles by his students and followers.
Attempts at formalizing the heuristic
ideas by Kantorovich  started at the initial stages
of $K$-space theory and yielded the so-called
identity preservation theorems.
They asserted that if some algebraic proposition with
finitely many function variables is satisfied by  the assignment of
all real values then
it remains valid  after replacement of reals with members of an~
arbitrary $K$-space.
To explain the nature of the Kantorovich transfer principle
became possible only after half a century by using
the technique of nonstandard models of set theory.

The abstract ideas of L.~V. Kantorovich in the area of $K$-spaces
are tied with linear programming and approximate methods of analysis.
He wrote about the still-unrevealed possibilities and underestimation
 of his theory for economics and remarked:

\hangindent 20pt{\small\quad\quad
But the comparison and correspondence relations
play an extraordinary role in economics and it was
definitely clear even at the cradle of $K$-spaces
that they will find their place in economic analysis and
yield luscious fruits.}
\medskip

The problem of  the scope of  the Hahn--Banach Theorem,
tantamount to describing the possible extensions of
linear programming, was rather popular in the decade past mid-1970s.
Everyone knows that linear programs lose their effectiveness
if only integer solutions are sought.
S.~N. Chernikov abstracted linear programming from the reals to some rings
similar to the rationals.\footnote{Cp.~\cite{Chernikov}.}
Rather topical in the world mathematical literature
was the problem of finding the algebraic systems that admit the full strength
of the ideas of L.~V. Kantorovich. The appropriate answer was given
by describing the abstract
modules that allow for the tools equivalent to the Hahn--Banach
Theorem.\footnote{Cp.~\cite{Modules}.}
These are $K$-spaces viewed as modules over rather ``voluminous''
algebras of their orthomorphisms. This result was resonated to some
extent in the theoretical background of mathematical economics as relevant to
the hypothesis of ``divisible goods.''

One of the rather simple particular cases of these results
is a theorem characterizing a lattice homomorphism. The latter
miraculously attracted attention of  vector-lattice theorists who
founded new proofs   and included the theorem in monographs as
{\it Kutateladze's Theorem}.\footnote{Cp.~\cite[p.~114]{Aliprantis}.}
Many years had elapsed before the usage of  Boolean valued models explained
that the modules found are in fact  dense subfields of the reals
in an appropriate nonstandard model of set theory.

In this area some unexpected generalizations of the Kre\u\i{}n--Milman
Theorem to noncompact sets that stimulate  a few articles
on the abstraction of Choquet theory to vector lattices.\footnote{Cp.~\cite{Choquet}.}
Ordered vector spaces have opened opportunities
to advance applications of Choquet theory to  several problems
of  modern potential theory such as describing  interconnections of the Dirichlet
problem with Bauer's geometric simplices in infinite dimensions
and introducing the new objects---supremal generators of function spaces 
which are convenient in approximation by positive operators.
Note that the conception of supremal generation
which bases on the computational simplicity of calculating the join of
two reals had turned out close to some ideas of idempotent analysis that
emerged somewhat later in the research by V.~P. Maslov and
his students.\footnote{In particular, cp.~\cite{Maslov}.}

\section{Nonsmooth Analysis and Optimization}

 Note the rather numerous papers on convex
analysis, one of the basic sections of nonlinear analysis.
Convex analysis is the calculus of linear inequalities.
The concept of convex set does not reach the age of 150 years, and convex
analysis as a branch of mathematics exists a bit longer than half a century.
The solution sets of simultaneous linear inequalities are the
same as convex sets which can be characterized by their gauges,
support functions or distributions of curvature.
Functional analysis is impossible without
convexity since the existence of a nonzero continuous linear functional
is provided if and only if the ambient space has nonempty proper
open convex subsets.

Convex surfaces have rather simple contingencies, and
convex functions are directionally differentiable in the natural sense
and their derivatives are nonlinear usually  in quite a few points.
But these points extreme in the direct and indirect senses are most
important. Study of the local behavior of possible fractures at extreme
points is the subject of subdifferential calculus.

The most general and complete formulas were found for
recalculating the values and solutions of  rather general convex extremal
problems under the changes of variables that preserve convexity.
The key to these formulas is the new trick of presenting an arbitrary convex
operator as the result of an affine change of variables in a particular
sublinear operator,  a member of some family enumerated by
cardinals.\footnote{The basic results in this area
were published in~\cite{Operators}. The literature uses the term
{\it Kutateladze's canonical operator} (cp.~\cite[pp.~123--125]{Rubinov} and
\cite[p.~92]{Tikhomirov}).}
These formulas led to the Lagrange principle for new classes of vector optimization
problems  and the theory of convex $\varepsilon$-programming.
The problem of approximate programming consists in the searching of a point
at which the value of a (possibly vector valued) function
differs from the extremum by at most  some positive error
vector $\varepsilon$. The constraints are also given to within
some accuracy of the order of $\varepsilon$.
The standard differential calculus is inapplicable here, but the new
methods of subdifferential calculus solve many problems of the sort.
These results became rather topical,
entered textbooks, and were redemonstrated with reference to the
Russian priority.\footnote{The literature uses the term
{\it Kutateladze's approximate solutions} (for instance, cp.~\cite{Novo}).}
Many years later the help of infinitesimal analysis
made it possible to propose the tricks that are not connected with
the bulky recalculations of errors. To this end, the error should be
considered as an infinitesimal, which is impossible within the classical
set-theoretic stance.

Applications to nonsmooth analysis are
connected primarily  with inspecting the
behavior of the contingencies of general rather than only convex correspondences.
In this area there were found some new rules for calculating
various types of tangents and one-sided directional derivatives.
The advances in these areas base on using the technique of model theory as well.

Many extremal problems are studied in
various branches of mathematics, but they use only scalar target functions.
Multiple criteria problems have appeared rather recently and
beyond the realm of mathematics. This explains the essential gap between
the complexity and efficiency of the mathematical tools 
which divides single and multiple criteria problems.
So it stands to reason to
enrich the stock of purely mathematical problems of vector optimization.
The author happened to distinguish  some class of geometrically reasonable
problems of vector optimization whose solutions can be presented in a
relatively lucid  form of conditions for surface area measures.
As model examples, the
Urysohn problems were considered  with extra targets like flattening in a given direction,
symmetry, or optimization of the volume of the convex hulls of several
surfaces.\footnote{Cp.~\cite{Multiple} and \cite{Balls}.}

\section{New Models for Mathematical Analysis}

   In the recent decades much research is done into the
 nonstandard methods located at the junctions
 of analysis and logic. This area requires the study of
some new opportunities of modeling that open broad vistas
 for consideration and solution of various theoretical and applied problems.

A model of a mathematical theory is usually called nonstandard if the membership
within the model has interpretation  different from that of set theory.\footnote{This understanding
is due to L.~Henkin.}  The simplest example of nonstandard modeling is the classical trick of
 presenting reals as points of  an axis.

The new methods of analysis are the adaptation of nonstandard set theoretic models
 to the problems of analysis. The two technologies are most popular: infinitesimal analysis
 also known as Robinsonian nonstandard analysis and Boolean valued analysis.

 Infinitesimal analysis by A.~Robinson appeared in 1960
and is characterized by legitimizing the usage of actual infinities and infinitesimals which were
forbidden for the span of about thirty years in the mathematics of the twentieth century.
In a sense, nonstandard analysis implements a partial modern return to the classical infinitesimal analysis.
The recent publications in this area can be partitioned into the two groups: The one that is most proliferous
uses infinitesimal analysis for ``killing quantifiers,'' i.e., simplifying definitions and proofs of the
classical results. The other has less instances but contributes much more to mathematics, searching
the opportunities unavailable to the standard methods; i.e., it develops the technologies
whose description is impossible without the new syntax based on the predicate of standardness.
We should list here the development of the new schemes for replacing 
the infinite objects as parts of finite sets: nonstandard hulls, Loeb measures, hyperapproximation, etc. 
Part of this research is done in Novosibirsk.
In particular, the author's results on infinitesimal programming belong to the second
group.\footnote{Cp.~\cite{InfProg}.}

Boolean valued analysis is characterized by the terms like Boolean valued universe,
descents and ascents, cyclic envelopes and mixings, Boolean sets and mappings, etc.
The technique here is much more complicated that of infinitesimal analysis
and just a few analysis are accustomed to it.
The rise of this branch of mathematical logic was connected with
the famous P.-J. Cohen's results of 1961 on the independence of the continuum hypothesis,
whose understanding drove D.~Scott, R. Solovay, and P.~Vop\v enka to the construction of the Boolean valued models of set theory.

D.~Scott foresaw the role of Boolean valued models in mathematics and wrote
   as far back as in~1969:\footnote{Cp.~\cite[p.~91]{Scott}.}

\begin{itemize}
\item[]{\small
We must ask whether there is
any interest in these nonstandard models aside from the independence proof;
that is, do they have any mathematical interest?
The answer must be yes, but we cannot yet give a really good argument.
}
\end{itemize}

G.~Takeuti
was one of the first who pointed out the role of these models
for functional analysis (in Hilbert space) and minted the term
{\it Boolean valued analysis}.\footnote{Cp.~\cite{Takeuti}.}
The models of infinitesimal analysis
can be viewed among the simplest instances of Boolean valued
universes.

The progress of Boolean valued analysis
in the recent decades has led to
a profusion of principally new ideas and results in many areas of
functional analysis, primary, in the theory of Dedekind complete
vector lattices and the theory of von Neumann algebras as well as
in convex analysis and the theory of vector measures.
Most of these advances are connected with
Novosibirsk.\footnote{Cp.~\cite{IBA}--\cite{INFA}.}
It is not an exaggeration to say that Boolean valued analysis left the
realm of logic and  has become a section of order analysis.

The new possibilities reveal the exceptional role
of universally complete vector lattices---extended $K$-spaces in the
Russian literature. It was completely unexpected that each of them turns out to be a legitimate
model of the real axis, so serving the same fundamental role in mathematics
as the reals.
Kantorovich spaces are indeed instances of the models
of the reals, which corroborated the heuristic ideas of
L.~V. Kantorovich.

Adaptation of nonstandard models to the problems of analysis
occupies the central place in the research of the author and his
closest colleagues. In this area we have developed the special
technique of  ascending and descending,  gave the criteria
of extensional algebraic systems, suggested the theory of cyclic
monads, and indicated some approaches to combining infinitesimal and
Boolean valued models.

These ideas lie behind solutions of various problems of geometric and applied
functional analysis among which we list
the drastically new classification of the Clarke type one-sided
approximations to arbitrary sets and the corresponding rules
for calculating infinitesimal tangents, the nonstandard approach
to approximate solutions of convex programs in the form of
{\it infinitesimal programming}, the new formulas for projecting to
the principal bands of the space of regular operators which are free
from the usual limitations on the order dual, etc.

We can also mention the new method of
studying some classes of bounded operators by the properties
of the kernels of their strata. This method bases on applying the
Kantorovich heuristic principle to the folklore fact  that a linear functional
can be restored from each of its hyperplanes to within a scalar multiplier.
In 2005 this led to the description of the operator annihilators of
Grothendieck  spaces.\footnote{Cp.~\cite{G-space}.}
In  2010 the method made it possible to suggest
the operator forms of the classical Farkas Lemma
in the theory of linear inequalities, so returning to the
origins of linear programming.\footnote{Cp.~\cite{Farkas} and \cite{Poly}.}

Of  great importance in this area are not only applications but also inspections
of the combined methods that involve Boolean valued and infinitesimal
techniques. At least the two approaches are viable:
One consists in studying a standard Boolean valued model
within the universes of Nelson's or Kawai's theory.
Infinitesimals descend there from
some external universe. The other bases on distinguishing infinitesimals within
Boolean valued models. These approaches were elaborated to some extent,
but the synthesis of the tools of various versions of nonstandard analysis
still remain a rather open problem.

Adaptation of the modern ideas of model theory
to functional analysis projects among the most important
directions of developing the synthetic methods of pure and applied
mathematics.  This approach yields new models of numbers,
spaces, and types of equations. The content expands of
all available theorems and algorithms. The whole methodology
of mathematical research is enriched and renewed, opening up
absolutely fantastic opportunities.
We can now use actual infinities and infinitesimals, transform
matrices into numbers, spaces into straight lines, and noncompact spaces into
compact spaces, yet having still uncharted vast territories of  new knowledge.

Quite a long time had passed until the classical functional analysis
occupied its present position of the language of continuous mathematics.
Now the time has come of the new powerful technologies of model theory
in mathematical analysis. Not all theoretical  and applied mathematicians
have already gained the importance of modern tools and
learned how to use them. However, there is no backward traffic in
science, and the modern methods
are doomed to reside in the realm of mathematics for ever and in a short
time they will become  as elementary and omnipresent in calculuses and
calculations as Banach spaces and linear operators.

\bibliographystyle{plain}

\begin{thebibliography}{99}

\bibitem{Reshetnyak}
Reshetnyak~Yu.~G.
{\it On the Length and Swerve of a~Curve  and the Area of a~Surface}
(PhD Thesis). Leningrad:  Leningrad State University, 1954.


\bibitem{Alexandrov}
Alexandrov~A.~D.
{\it Selected Works. Vol.~1: Geometry and Applications.}
Novosibirsk: Nauka, 2006.


\bibitem{Duality}
Kutateladze~S.~S. and Rubinov~A.~M.
``Minkowski duality and its applications.''
{\it Russian Math. Surveys},  {\bf27}:3, 137--191 (1972).


\bibitem{Singer}
Singer I.
{\it Abstract Convex Analysis.} New York: John Wiley, 1997.



\bibitem{Urysohn}
Urysohn~P.~S. (1924)
``Dependence between the average width and volume of convex bodies.''
{\it Mat. Sb.,} {\bf31}:3,  477--486 (1924).



\bibitem{Pogorelov}
Pogorelov~A.~V.
``Imbedding a `soap bubble' into a tetrahedron.''
{\it Math. Notes}, {\bf56}:2, 824--826 (1994).


\bibitem{LVK}
Kantorovich L.~V.
``On semiordered linear spaces and their
applications in the theory of linear operators,''
{\it Dokl. Akad. Nauk SSSR}, {\bf 4}: 1--2, 11--14 (1935).



\bibitem{Chernikov}
Chernikov S.~N.
{\it Linear Inequalities}. Moscow: Nauka, 1968.


\bibitem{Modules}
Kutateladze S.~S.
``Modules admitting convex analysis.''
{\it Soviet Math. Dokl.,} {\bf 21}:3, 820-823 (1980).


\bibitem{Aliprantis}
Aliprants Ch. and Birkinshaw Ow.
{\it Positive Operators.} Orlando etc.: Academic Press, 1985.


\bibitem{Choquet}
Kutateladze S.~S.
``Choquet boundaries in K-spaces.''
{\it Russian Math. Surveys.} 1975, {\bf30}:4, 115--155 (1975).


\bibitem{Maslov}
Kolokol'tsov V.~N. and Maslov V.~P.
``Idempotent analysis as a tool of control theory and optimal synthesis. I.''
{\it Funct. Anal. Appl.}, {\bf23}:1, 1--11 (1989).


\bibitem{Operators}
Kutateladze S.~S.
``Convex operators.''
{\it Russian Math. Surveys}, {\bf34}:1, 181-214 (1979).

\bibitem{Rubinov}
Rubinov A.~M.
``Sublinear operators and their applications.''
{\it Russian Math. Surveys}, {\bf32}:4, 115--175 (1977).

\bibitem{Tikhomirov}
Tikhomirov V. M.
``Convex analysis.''
In: {\it Current Problems in Mathematics. Fundamental Trends,} {\bf14}. Moscow: VINITI,
5--101, 1987.



\bibitem{Novo}
Guti\'errez C., Jim\'enez B., and Novo V.
``On approximate solutions in vector optimization
problems via scalarization.''
{\it Computat. Optim. Appl.}, {\bf35}, 305--324 (2006).

\bibitem{Multiple}
Kutateladze S.~S.
``Multiobjective problems of convex geometry.''
{\it Siberian Math. J.,} {\bf50}:5, 887--897 (2009.)



\bibitem{Balls}
Kutateladze S.~S.
``Multiple criteria problems over Minkowski balls.''
{\it J.~Appl.~Indust.~Math.,} {\bf7}:2, 208--214 (2013).

\bibitem{InfProg}
Kutateladze S.~S.
``A variant of nonstandard convex programming.''
{\it Siberian Math. J.,} {\bf27}: 4, 537--544 (1986).


\bibitem{Scott}
Scott D.
``Boolean models and nonstandard analysis.''
In: {\it Applications of Model Theory to Algebra, Analysis, and Probability}, 87--92.
New York: Holt, Rinehart and Winston (1969).


\bibitem{Takeuti}
Takeuti  G.
{\it Two Applications of Logic to Mathematics.}
Tokyo and Princeton: Iwanami Publ. \& Princeton University Press, 1978.



\bibitem{IBA}
Kusraev A.~G. and Kutateladze S.~S.
{\it Introduction to Booolean Valued Analysis.}
Moscow: Nauka, 2005.

\bibitem{BA_Selected}
Kusraev A.~G. and Kutateladze S.~S.
{\it Boolean Valued Analysis: Selected Topics.}
Vladikavkaz: Southern Math. Institute, 2014.



\bibitem{INFA}
Gordon E.~I., Kusraev A.~G., and Kutateladze S.~S.
{\it Infinitesimal Analysis: Selected Topics}.
Moscow: Nauka, 2011.



\bibitem{G-space}
Kutateladze S.~S.
``On Grothendieck subspaces.''
{\it Siberian Math. J.,}  {\bf46}:3, 489--493 (2005).



\bibitem{Farkas}
Kutateladze S.~S.
``The Farkas lemma revisited.''
{\it Siberian Math. J.,}
 2010, {\bf51}:1, 78--87 (2010).

\bibitem{Poly}
Kutateladze S.~S.
``The polyhedral Lagrange principle,''
{\it Siberian Math. J.} 2011, {\bf52}:3, 484--486 (2011).

\end{thebibliography}

\end{document}